\documentclass{amsart}
\usepackage{amsthm}
\usepackage{amstext}
\usepackage{amssymb}

\makeatletter

\numberwithin{equation}{section}
\numberwithin{figure}{section}
  \theoremstyle{remark}
  \newtheorem*{note*}{Note}
\theoremstyle{plain}
\newtheorem{thm}{Theorem}
\newcommand\relphantom[1]{\mathrel{\phantom{{#1}}}}

\makeatother

\allowdisplaybreaks
\usepackage[sort&compress,square,numbers]{natbib}

\begin{document}

\title{Fully degenerate poly-Bernoulli numbers and polynomials}
\author{Dae San Kim}
\address{Department of Mathematics, Sogang University, Seoul 121-742, Republic
of Korea}
\email{dskim@sogang.ac.kr}

\author{Taekyun Kim}
\address{Department of Mathematics, Kwangwoon University, Seoul 139-701, Republic
of Korea}
\email{tkkim@kw.ac.kr}


\begin{abstract}
In this paper, we introduce the new fully degenerate poly-Bernoulli numbers
and polynomials and inverstigate some properties of these polynomials
and numbers. From our properties, we derive some identities for the
fully degenerate poly-Bernoulli numbers and polynomials.
\end{abstract}

\maketitle

\global\long\def\Zp{\mathbb{Z}_{p}}

\global\long\def\li{\mathrm{Li}}

\global\long\def\acl#1#2{\left\langle \left.#1\right|#2\right\rangle }

\global\long\def\acr#1#2{\left\langle #1\left|#2\right.\right\rangle }

\section{Introduction}

It is well known that the Bernoulli polynomials are defined by the
generating function
\begin{equation}
\frac{t}{e^{t}-1}e^{xt}=\sum_{n=0}^{\infty}B_{n}\left(x\right)\frac{t^{n}}{n!},\quad\left(\text{see
\cite{key-1,key-2,key-3,key-4,key-5,key-6,key-7,key-8,key-9,key-10,key-11,key-12,key-13,key-14,key-15,key-16,key-17,key-18,key-19,key-20,key-21}}\right).\label{eq:1}
\end{equation}
When $x=0$, $B_{n}=B_{n}\left(0\right)$ are called the Bernoulli
numbers. From (\ref{eq:1}), we note that
\begin{equation}
B_{n}\left(x\right)=\sum_{l=0}^{n}\binom{n}{l}B_{l}x^{n-l},\quad\left(n\ge0\right),\label{eq:2}
\end{equation}
and
\begin{equation}
B_{0}=1,\quad B_{n}\left(1\right)-B_{n}=\delta_{1,n},\quad\left(n\in\mathbb{N}\right),\quad\left(\text{see \cite{key-1,key-19}}\right),\label{eq:3}
\end{equation}
where $\delta_{n,k}$ is the Kronecker's symbol.

In \cite{key-3}, L. Carlitz considered the degenerate Bernoulli polynomials
which are given by the generating function
\begin{equation}
\frac{t}{\left(1+\lambda t\right)^{\frac{1}{\lambda}}-1}\left(1+\lambda t\right)^{\frac{x}{\lambda}}=\sum_{n=0}^{\infty}\beta_{n,\lambda}\left(x\right)\frac{t^{n}}{n!}.\label{eq:4}
\end{equation}

When $x=0$, $\beta_{n,\lambda}=\beta_{n,\lambda}\left(0\right)$
are called the degenerate Bernoulli numbers. From (\ref{eq:1}) and
(\ref{eq:4}), we note that
\begin{equation}
\frac{t}{e^{t}-1}e^{xt}=\lim_{\lambda\rightarrow0}\frac{t}{\left(1+\lambda t\right)^{\frac{1}{\lambda}}-1}\left(1+\lambda t\right)^{\frac{x}{\lambda}}=\sum_{n=0}^{\infty}\lim_{\lambda\rightarrow0}\beta_{n,\lambda}\left(x\right)\frac{t^{n}}{n!.}\label{eq:5}
\end{equation}

Thus, by (\ref{eq:5}), we get
\begin{equation}
B_{n}\left(x\right)=\lim_{\lambda\rightarrow0}\beta_{n,\lambda}\left(x\right),\quad\left(\text{see \cite{key-3,key-14}}\right).\label{eq:6}
\end{equation}

By (\ref{eq:4}), we get
\begin{equation}
\sum_{n=0}^{\infty}\left(\beta_{m,\lambda}\left(1\right)-\beta_{m,\lambda}\right)\frac{t^{n}}{n!}=t,\label{eq:7}
\end{equation}
and
\begin{equation}
\beta_{n,\lambda}\left(x\right)=\sum_{l=0}^{n}\binom{n}{l}\beta_{l,\lambda}\lambda^{n-l}\left(\frac{x}{\lambda}\right)_{n-l}.\label{eq:8}
\end{equation}

From (\ref{eq:7}), we have
\begin{equation}
\beta_{n,\lambda}\left(1\right)-\beta_{n,\lambda}=\delta_{1,n},\quad\left(n\ge0\right),\quad\beta_{0,\lambda}=1.\label{eq:9}
\end{equation}

Now, we consider the degenerate Bernoulli polynomials which are different from
the degenerate Bernoulli polynomials of L. Carlitz as follows:
\begin{equation}
\frac{\log\left(1+\lambda t\right)^{\frac{1}{\lambda}}}{\left(1+\lambda t\right)^{\frac{1}{\lambda}}-1}\left(1+\lambda t\right)^{\frac{x}{\lambda}}=\sum_{n=0}^{\infty}b_{n,\lambda}\left(x\right)\frac{t^{n}}{n!}.\label{eq:10}
\end{equation}

When $x=0$, $b_{n,\lambda}=b_{n,\lambda}\left(0\right)$ are called
the degenerate Bernoulli numbers.
\begin{note*}
The degenerate Bernoulli polynomials are also called Daehee
polynomials with $\lambda$-parameter (see \cite{key-13}).
\end{note*}
From (\ref{eq:10}), we note that
\begin{align}
\frac{t}{e^{t}-1}e^{xt} & =\lim_{\lambda\rightarrow0}\frac{\log\left(1+\lambda t\right)^{\frac{1}{\lambda}}}{\left(1+\lambda t\right)^{\frac{1}{\lambda}}-1}\left(1+\lambda t\right)^{\frac{x}{\lambda}}\label{eq:11}\\
 & =\sum_{n=0}^{\infty}\lim_{\lambda\rightarrow0}b_{n,\lambda}\left(x\right)\frac{t^{n}}{n!}.\nonumber
\end{align}

By (\ref{eq:1}) and (\ref{eq:11}), we see that
\[
B_{n}\left(x\right)=\lim_{\lambda\rightarrow0}b_{n,\lambda}\left(x\right),\quad\left(n\ge0\right).
\]

The classical polylogarithm function $\li_{k}\left(x\right)$ is defined
by
\begin{equation}
\li_{k}\left(t\right)=\sum_{n=1}^{\infty}\frac{t^{n}}{n^{k}},\quad\left(k\in\mathbb{Z}\right),\quad\left(\text{see \cite{key-10,key-11}}\right).\label{eq:12}
\end{equation}

It is known that the poly-Bernoulli polynomials are defined by the
generating function
\begin{equation}
\frac{\li_{k}\left(1-e^{-t}\right)}{1-e^{-t}}e^{xt}=\sum_{n=0}^{\infty}B_{n}^{\left(k\right)}\left(x\right)\frac{t^{n}}{n!},\quad\left(\text{see \cite{key-9,key-10,key-16}}\right).\label{eq:13}
\end{equation}

When $k=1$, we have
\begin{align}
\sum_{n=0}^{\infty}B_{n}^{\left(1\right)}\left(x\right)\frac{t^{n}}{n!} & =\frac{t}{1-e^{-t}}e^{xt}=\frac{t}{e^{t}-1}e^{\left(x+1\right)t}\label{eq:14}\\
 & =\sum_{n=0}^{\infty}B_{n}\left(x+1\right)\frac{t^{n}}{n!}.\nonumber
\end{align}

By (\ref{eq:14}), we easily get
\[
B_{n}^{\left(1\right)}\left(x\right)=B_{n}\left(x+1\right),\quad\left(n\ge0\right).
\]

Let $x=0$. Then $B_{n}^{\left(k\right)}=B_{n}^{\left(k\right)}\left(0\right)$
are called the poly-Bernoulli numbers.

In this paper, we introduce the new fully degenerate poly-Bernoulli numbers
and polynomials and inverstigate some properties of these polynomials
and numbers. From our investigation, we derive some identities for
the fully degenerate poly-Bernoulli numbers and polynomials.

\section{Fully degenerate poly-Bernoulli polynomials}

For $k\in\mathbb{Z}$, we define the fully degenerate poly-Bernoulli polynomials
which are given by the generating function
\begin{equation}
\frac{\li_{k}\left(1-\left(1+\lambda t\right)^{-\frac{1}{\lambda}}\right)}{1-\left(1+\lambda t\right)^{-\frac{1}{\lambda}}}\left(1+\lambda t\right)^{\frac{x}{\lambda}}=\sum_{n=0}^{\infty}\beta_{n,\lambda}^{\left(k\right)}\left(x\right)\frac{t^{n}}{n!}.\label{eq:15}
\end{equation}

When $x=0$, $\beta_{n,\lambda}^{\left(k\right)}=\beta_{n,\lambda}^{\left(k\right)}\left(0\right)$
are called the fully degenerate poly-Bernoulli numbers.

From (\ref{eq:13}) and (\ref{eq:15}), we have
\begin{align}
\frac{\li_{k}\left(1-e^{-t}\right)}{1-e^{-t}}e^{xt} & =\lim_{\lambda\rightarrow0}\frac{\li_{k}\left(1-\left(1+\lambda t\right)^{-\frac{1}{\lambda}}\right)}{1-\left(1+\lambda t\right)^{-\frac{1}{\lambda}}}\left(1+\lambda t\right)^{\frac{x}{\lambda}}\label{eq:16}\\
 & =\sum_{n=0}^{\infty}\lim_{\lambda\rightarrow0}\beta_{n,\lambda}^{\left(k\right)}\left(x\right)\frac{t^{n}}{n!}.\nonumber
\end{align}

Thus, we get
\begin{equation}
\lim_{\lambda\rightarrow0}\beta_{n,\lambda}^{\left(k\right)}\left(x\right)=B_{n}^{\left(k\right)}\left(x\right),\quad\left(n\ge0\right).\label{eq:17}
\end{equation}

By (\ref{eq:15}), we get
\begin{align}
 & \relphantom{=}\sum_{n=0}^{\infty}\beta_{n,\lambda}^{\left(k\right)}\left(x\right)\frac{t^{n}}{n!}\label{eq:18}\\
 & =\frac{\li_{k}\left(1-\left(1+\lambda t\right)^{-\frac{1}{\lambda}}\right)}{1-\left(1+\lambda t\right)^{-\frac{1}{\lambda}}}\left(1+\lambda t\right)^{\frac{x}{\lambda}}\nonumber \\
 & =\sum_{n=0}^{\infty}\left(\sum_{l=0}^{n}\binom{n}{l}\beta_{l,\lambda}^{\left(k\right)}\left(\frac{x}{\lambda}\right)_{n-l}\lambda^{n-l}\right)\frac{t^{n}}{n!}.\nonumber \\
\nonumber
\end{align}

Thus, from (\ref{eq:18}), we have
\begin{equation}
\beta_{n,\lambda}^{\left(k\right)}\left(x+y\right)=\sum_{l=0}^{n}\binom{n}{l}\left(\frac{y}{\lambda}\right)_{n-l}\lambda^{n-l}\beta_{l,\lambda}^{\left(k\right)}\left(x\right),\quad\left(n\ge0\right),\label{eq:19}
\end{equation}
and
\[
\beta_{n,\lambda}^{\left(k\right)}\left(x\right)=\sum_{l=0}^{n}\binom{n}{l}\left(\frac{x}{\lambda}\right)_{n-l}\lambda^{n-l}\beta_{l,\lambda}^{\left(k\right)}.
\]

Therefore, by (\ref{eq:17}) and (\ref{eq:19}), we obtain the following
theorem.
\begin{thm}
\label{thm:1} For $k\in\mathbb{Z}$, $n\ge0$, we have
\begin{equation}
\beta_{n,\lambda}^{\left(k\right)}\left(x+y\right)=\sum_{l=0}^{n}\binom{n}{l}\left(\frac{y}{\lambda}\right)_{n-l}\lambda^{n-l}\beta_{l,\lambda}^{\left(k\right)}\left(x\right),\quad\left(n\ge0\right),\label{eq:19-1}
\end{equation}
and
\[
\lim_{\lambda\rightarrow0}\beta_{n,\lambda}^{\left(k\right)}\left(x\right)=B_{n}^{\left(k\right)}\left(x\right),
\]
where $\left(x\right)_{n}=x\left(x-1\right)\cdots\left(x-n+1\right)=\sum_{l=0}^{n}S_{1}\left(n,l\right)x^{l}$.
\end{thm}
From (\ref{eq:15}), we can derive the following equation:
\begin{equation}
\sum_{n=0}^{\infty}\beta_{n,\lambda}^{\left(k\right)}\left(x\right)\frac{t^{n}}{n!}=\frac{\li_{k}\left(1-\left(1+\lambda t\right)^{-\frac{1}{\lambda}}\right)}{\left(1+\lambda t\right)^{\frac{1}{\lambda}}-1}\left(1+\lambda t\right)^{\frac{x+1}{\lambda}}.\label{eq:20}
\end{equation}

Thus, by (\ref{eq:20}), we get
\begin{align}
 & \relphantom{=}\sum_{n=0}^{\infty}\left\{ \beta_{n,\lambda}^{\left(k\right)}-\beta_{n,\lambda}^{\left(k\right)}\left(-1\right)\right\} \frac{t^{n}}{n!}\label{eq:21}\\
 & =\li_{k}\left(1-\left(1+\lambda t\right)^{-\frac{1}{\lambda}}\right)\nonumber \\
 & =\sum_{m=1}^{\infty}\frac{\left(1-\left(1+\lambda t\right)^{-\frac{1}{\lambda}}\right)^{m}}{m^{k}}\nonumber \\
 & =\sum_{m=0}^{\infty}\frac{\left(-1\right)^{m+1}}{\left(m+1\right)^{k}}\left(e^{-\frac{1}{\lambda}\log\left(1+\lambda t\right)}-1\right)^{m+1}\nonumber \\
 & =\sum_{m=0}^{\infty}\frac{\left(-1\right)^{m+1}}{\left(m+1\right)^{k}}\left(m+1\right)!\sum_{l=m+1}^{\infty}S_{2}\left(l,m+1\right)\left(-1\right)^{l}\lambda^{-l}\frac{\left(\log\left(1+\lambda t\right)\right)^{l}}{l!}\nonumber \\
 & =\sum_{l=1}^{\infty}\sum_{m=0}^{l-1}\frac{\left(-1\right)^{m+1}}{\left(m+1\right)^{k}}\left(m+1\right)!S_{2}\left(l,m+1\right)\left(-1\right)^{l}\lambda^{-l}\sum_{n=l}^{\infty}S_{1}\left(n,l\right)\lambda^{n}\frac{t^{n}}{n!}\nonumber \\
 & =\sum_{n=1}^{\infty}\left(\sum_{l=1}^{n}\sum_{m=0}^{l-1}\frac{m!\left(m+1\right)\left(-1\right)^{l-m-1}\lambda^{n-l}S_{2}\left(l,m+1\right)S_{1}\left(n,l\right)}{\left(m+1\right)^{k}}\right)\frac{t^{n}}{n!},\nonumber
\end{align}
where $S_{2}\left(n,l\right)$ and $S_{1}\left(n,l\right)$ are the Stirling
numbers of the second kind and of the first kind, respectively.

Therefore, by (\ref{eq:21}), we obtain the following theorem.
\begin{thm}
\label{thm:2} For $k\in\mathbb{Z}$, $n\ge1$, we have
\[
\beta_{n,\lambda}^{\left(k\right)}-\beta_{n,\lambda}^{\left(k\right)}\left(-1\right)=\sum_{l=1}^{n}\sum_{m=0}^{l-1}\frac{m!\left(-1\right)^{l-m-1}\lambda^{n-l}S_{2}\left(l,m+1\right)S_{1}\left(n,l\right)}{\left(m+1\right)^{k-1}}.
\]

\end{thm}
From (\ref{eq:12}), we can easily derive the following equation:
\begin{equation}
\li_{k}^{\prime}\left(t\right)=\frac{d}{dt}\li_{k}\left(t\right)=\frac{1}{t}\li_{k-1}\left(t\right).\label{eq:22}
\end{equation}

Thus, by (\ref{eq:22}), the generating function of the fully degenerate
poly-Bernoulli numbers is also written in terms of the following iterated integral:
\begin{align}
 & \frac{\left(1+\lambda t\right)^{\frac{1}{\lambda}}}{\left(1+\lambda t\right)^{\frac{1}{\lambda}}-1}\int_{0}^{t}\frac{1}{\left(\left(1+\lambda t\right)^{\frac{1}{\lambda}}-1\right)\left(1+\lambda t\right)}\label{eq:23}\\
 &\times \int_{0}^{t}\frac{1}{\left(\left(1+\lambda t\right)^{\frac{1}{\lambda}}-1\right)\left(1+\lambda t\right)}\cdots\int_{0}^{t}\frac{\log\left(1+\lambda t\right)^{\frac{1}{\lambda}}}{\left(\left(1+\lambda t\right)^{\frac{1}{\lambda}}-1\right)\left(1+\lambda t\right)}\underset{k-1\text{ times}}{\underbrace{dt\cdots dt}}\nonumber\\
 & =\sum_{n=0}^{\infty}\beta_{n,\lambda}^{\left(k\right)}\frac{t^{n}}{n!}.\nonumber
\end{align}

For $k=2$, we have
\begin{align}
 & \relphantom{=}\sum_{n=0}^{\infty}\beta_{n,\lambda}^{\left(2\right)}\frac{t^{n}}{n!}\label{eq:24}\\
 & =\frac{\left(1+\lambda t\right)^{\frac{1}{\lambda}}}{\left(1+\lambda t\right)^{\frac{1}{\lambda}}-1}\int_{0}^{t}\frac{\log\left(1+\lambda t\right)^{\frac{1}{\lambda}}}{\left(1+\lambda t\right)^{\frac{1}{\lambda}}-1}\left(1+\lambda t\right)^{-\frac{\lambda}{\lambda}}dt\nonumber \\
 & =\frac{\left(1+\lambda t\right)^{\frac{1}{\lambda}}}{\left(1+\lambda t\right)^{\frac{1}{\lambda}}-1}\left(\sum_{m=0}^{\infty}b_{m,\lambda}\left(-\lambda\right)\frac{1}{m!}\int_{0}^{t}t^{m}dt\right)\nonumber \\
 & =\left(\frac{t}{\left(1+\lambda t\right)^{\frac{1}{\lambda}}-1}\left(1+\lambda t\right)^{\frac{1}{\lambda}}\right)\left(\sum_{m=0}^{\infty}\frac{b_{m,\lambda}\left(-\lambda\right)}{\left(m+1\right)}\frac{t^{m}}{m!}\right)\nonumber \\
 & =\sum_{n=0}^{\infty}\left(\sum_{l=0}^{n}\binom{n}{l}\beta_{l,\lambda}\left(1\right)\frac{b_{n-l,\lambda}\left(-\lambda\right)}{n-l+1}\right)\frac{t^{n}}{n!}.\nonumber
\end{align}

Therefore, by (\ref{eq:24}), we obtain the following theorem.
\begin{thm}
\label{thm:3} For $n\ge0$, we have
\[
\beta_{n,\lambda}^{\left(2\right)}=\sum_{l=0}^{n}\binom{n}{l}\beta_{l,\lambda}\left(1\right)\frac{b_{n-l,\lambda}\left(-\lambda\right)}{n-l+1}.
\]

\end{thm}
Note that
\[
B_{n}^{\left(2\right)}=\lim_{\lambda\rightarrow0}\beta_{n,\lambda}^{\left(2\right)}=\sum_{l=0}^{n}\binom{n}{l}B_{l}\left(1\right)\frac{B_{n-l}}{n-l+1}.
\]

From (\ref{eq:15}), we have
\begin{align}
 & \relphantom{=}\sum_{n=0}^{\infty}\beta_{n,\lambda}^{\left(k\right)}\frac{t^{n}}{n!}\label{eq:25}\\
 & =\frac{\li_{k}\left(1-\left(1+\lambda t\right)^{-\frac{1}{\lambda}}\right)}{1-\left(1+\lambda t\right)^{-\frac{1}{\lambda}}}\nonumber \\
 & =\sum_{m=0}^{\infty}\frac{1}{\left(m+1\right)^{k}}\left(1-\left(1+\lambda t\right)^{-\frac{1}{\lambda}}\right)^{m}\nonumber \\
 & =\sum_{m=0}^{\infty}\frac{\left(-1\right)^{m}}{\left(m+1\right)^{k}}\left(e^{-\frac{1}{\lambda}\log\left(1+\lambda t\right)}-1\right)^{m}\nonumber \\
 & =\sum_{m=0}^{\infty}\frac{\left(-1\right)^{m}}{\left(m+1\right)^{k}}m!\sum_{l=m}^{\infty}S_{2}\left(l,m\right)\left(-\frac{1}{\lambda}\right)^{l}\frac{\left(\log\left(1+\lambda t\right)\right)^{l}}{l!}\nonumber \\
 & =\sum_{l=0}^{\infty}\left(\sum_{m=0}^{l}\frac{\left(-1\right)^{m+l}m!}{\left(m+1\right)^{k}}S_{2}\left(l,m\right)\lambda^{-l}\right)\frac{1}{l!}\left(\log\left(1+\lambda t\right)\right)^{l}\nonumber \\
 & =\sum_{l=0}^{\infty}\left(\sum_{m=0}^{l}\frac{\left(-1\right)^{m+l}m!}{\left(m+1\right)^{k}}S_{2}\left(l,m\right)\lambda^{-l}\right)\sum_{n=l}^{\infty}S_{1}\left(n,l\right)\lambda^{n}\frac{t^{n}}{n!}\nonumber \\
 & =\sum_{n=0}^{\infty}\left\{ \sum_{l=0}^{n}\sum_{m=0}^{l}\frac{\left(-1\right)^{m+l}m!}{\left(m+1\right)^{k}}S_{2}\left(l,m\right)S_{1}\left(n,l\right)\lambda^{n-l}\right\} \frac{t^{n}}{n!}.\nonumber
\end{align}

Therefore, by (\ref{eq:25}), we obtain the following theorem.
\begin{thm}
\label{thm:4} For $n\ge0$, we have
\[
\beta_{n,\lambda}^{\left(k\right)}=\sum_{l=0}^{n}\sum_{m=0}^{l}\frac{\left(-1\right)^{m+l}m!}{\left(m+1\right)^{k}}S_{2}\left(l,m\right)S_{1}\left(n,l\right)\lambda^{n-l}.
\]

\end{thm}
Note that
\[
B_{n}^{\left(k\right)}=\lim_{\lambda\rightarrow0}\beta_{n,\lambda}^{\left(k\right)}=\sum_{m=0}^{n}\frac{\left(-1\right)^{m+n}m!}{\left(m+1\right)^{k}}S_{2}\left(n,m\right).
\]

From (\ref{eq:22}), we have
\begin{align}
 & \relphantom{=}\frac{d}{dt}\li_{k}\left(1-\left(1+\lambda t\right)^{-\frac{1}{\lambda}}\right)\label{eq:28}\\
 & =\frac{1}{1-\left(1+\lambda t\right)^{-\frac{1}{\lambda}}}\left(1+\lambda t\right)^{-\frac{1}{\lambda}-1}\li_{k-1}\left(1-\left(1+\lambda t\right)^{-\frac{1}{\lambda}}\right)\nonumber \\
 & =\left(1+\lambda t\right)^{-\frac{1}{\lambda}-1}\sum_{n=0}^{\infty}\beta_{n,\lambda}^{\left(k-1\right)}\frac{t^{n}}{n!}.\nonumber
\end{align}

On the other hand,
\begin{align}
 & \relphantom{=}\frac{d}{dt}\left(\li_{k}\left(1-\left(1+\lambda t\right)^{-\frac{1}{\lambda}}\right)\right)\label{eq:29}\\
 & =\frac{d}{dt}\left(\left(1-\left(1+\lambda t\right)^{-\frac{1}{\lambda}}\right)\frac{1}{1-\left(1+\lambda t\right)^{-\frac{1}{\lambda}}}\li_{k}\left(1-\left(1+\lambda t\right)^{-\frac{1}{\lambda}}\right)\right)\nonumber \\
 & =\left(1+\lambda t\right)^{-\frac{1}{\lambda}-1}\frac{1}{1-\left(1+\lambda t\right)^{-\frac{1}{\lambda}}}\li_{k}\left(1-\left(1+\lambda t\right)^{-\frac{1}{\lambda}}\right)\nonumber \\
 & \relphantom{=}+\left(1-\left(1+\lambda t\right)^{-\frac{1}{\lambda}}\right)\frac{d}{dt}\left(\sum_{n=0}^{\infty}\beta_{n,\lambda}^{\left(k\right)}\frac{t^{n}}{n!}\right)\nonumber \\
 & =\left(1+\lambda t\right)^{-\frac{1}{\lambda}-1}\sum_{n=0}^{\infty}\beta_{n,\lambda}^{\left(k\right)}\frac{t^{n}}{n!}\nonumber \\
 & \relphantom{=}+\left(1-\left(1+\lambda t\right)^{-\frac{1}{\lambda}}\right)\sum_{n=1}^{\infty}\beta_{n,\lambda}^{\left(k\right)}\frac{t^{n-1}}{\left(n-1\right)!}.\nonumber
\end{align}

By (\ref{eq:28}) and (\ref{eq:29}), we get
\begin{align}
 & \relphantom{=}\sum_{n=0}^{\infty}\beta_{n,\lambda}^{\left(k-1\right)}\frac{t^{n}}{n!}\label{eq:30}\\
 & =\sum_{n=0}^{\infty}\beta_{n,\lambda}^{\left(k\right)}\frac{t^{n}}{n!}+\left(1+\lambda t\right)\left(\left(1+\lambda t\right)^{\frac{1}{\lambda}}-1\right)\sum_{n=1}^{\infty}\beta_{n,\lambda}^{\left(k\right)}\frac{t^{n-1}}{\left(n-1\right)!}\nonumber \\
 & =\sum_{n=0}^{\infty}\beta_{n,\lambda}^{\left(k\right)}\frac{t^{n}}{n!}+\left(\left(1+\lambda t\right)^{\frac{1}{\lambda}}-1\right)\sum_{m=0}^{\infty}\beta_{m+1,\lambda}^{\left(k\right)}\frac{t^{m}}{m!}\nonumber \\
 & \relphantom{=}+\lambda\left(\left(1+\lambda t\right)^{\frac{1}{\lambda}}-1\right)\sum_{m=0}^{\infty}\beta_{m,\lambda}^{\left(k\right)}m\frac{t^{m}}{m!}\nonumber \\
 & =\sum_{n=0}^{\infty}\beta_{n,\lambda}^{\left(k\right)}\frac{t^{n}}{n!}+\left(\sum_{l=1}^{\infty}\left(\frac{1}{\lambda}\right)_{l}\lambda^{l}\frac{t^{l}}{l!}\right)\left(\sum_{m=0}^{\infty}\beta_{m+1,\lambda}^{\left(k\right)}\frac{t^{m}}{m!}\right)\nonumber \\
 & \relphantom{=}+\lambda\left(\sum_{l=1}^{\infty}\left(\frac{1}{\lambda}\right)_{l}\lambda^{l}\frac{t^{l}}{l!}\right)\left(\sum_{m=0}^{\infty}\beta_{m,\lambda}^{\left(k\right)}m\frac{t^{m}}{m!}\right)\nonumber \\
 & =\sum_{n=0}^{\infty}\beta_{n,\lambda}^{\left(k\right)}\frac{t^{n}}{n!}+\sum_{n=1}^{\infty}\left(\sum_{m=0}^{n-1}\left(\frac{1}{\lambda}\right)_{n-m}\lambda^{n-m}\frac{n!}{\left(n-m\right)!m!}\beta_{m+1,\lambda}^{\left(k\right)}\right)\frac{t^{n}}{n!}\nonumber \\
 & \relphantom{=}+\lambda\sum_{n=1}^{\infty}\left(\sum_{m=0}^{n-1}\left(\frac{1}{\lambda}\right)_{n-m}\lambda^{n-m}\frac{m\cdot n!}{\left(n-m\right)!m!}\beta_{m,\lambda}^{\left(k\right)}\right)\frac{t^{n}}{n!}\nonumber \\
 & =\sum_{n=0}^{\infty}\beta_{n,\lambda}^{\left(k\right)}\frac{t^{n}}{n!}+\sum_{n=1}^{\infty}\left(\sum_{m=0}^{n-1}\binom{n}{m}\left(\frac{1}{\lambda}\right)_{n-m}\lambda^{n-m}\beta_{m+1,\lambda}^{\left(k\right)}\right)\frac{t^{n}}{n!}\nonumber \\
 & \relphantom{=}+\lambda\sum_{n=1}^{\infty}\left(\sum_{m=0}^{n-1}\binom{n}{m}\left(\frac{1}{\lambda}\right)_{n-m}\lambda^{n-m}m\beta_{m,\lambda}^{\left(k\right)}\right)\frac{t^{n}}{n!},\nonumber \\
\nonumber
\end{align}
where
\[
\left(\frac{1}{\lambda}\right)_{n}=\left(\frac{1}{\lambda}\right)\left(\frac{1}{\lambda}-1\right)\cdots\left(\frac{1}{\lambda}-n+1\right)=\sum_{l=0}^{n}S_{1}\left(n,l\right)\lambda^{-l},\quad\left(n\ge0\right).
\]

Thus, by (\ref{eq:30}), we have
\begin{align}
\beta_{n,\lambda}^{\left(k-1\right)} & =\beta_{n,\lambda}^{\left(k\right)}+\sum_{m=0}^{n-1}\binom{n}{m}\left(\frac{1}{\lambda}\right)_{n-m}\lambda^{n-m}\beta_{m+1,\lambda}^{\left(k\right)}\label{eq:31}\\
 & \relphantom{=}+\lambda\sum_{m=0}^{n-1}\binom{n}{m}\left(\frac{1}{\lambda}\right)_{n-m}\lambda^{n-m}m\beta_{m,\lambda}^{\left(k\right)}\nonumber \\
 & =\left(n+1\right)\beta_{n,\lambda}^{\left(k\right)}+\sum_{m=1}^{n-1}\binom{n}{m-1}\left(\frac{1}{\lambda}\right)_{n-m+1}\lambda^{n-m+1}\beta_{m,\lambda}^{\left(k\right)}\nonumber \\
 & \relphantom{=}+\lambda\sum_{m=0}^{n-1}\binom{n}{m}\left(\frac{1}{\lambda}\right)_{n-m}\lambda^{n-m}m\beta_{m,\lambda}^{\left(k\right)},\quad\left(n\ge1\right).\nonumber
\end{align}

Therefore, by (\ref{eq:31}), we obtain the following theorem.
\begin{thm}
\label{thm:5} For $n\ge1$, we have
\begin{align*}
\beta_{n,\lambda}^{\left(k\right)}&=\frac{1}{n+1}\left\{ \beta_{n,\lambda}^{\left(k-1\right)}-\sum_{m=1}^{n-1}\binom{n}{m-1}\beta_{m,\lambda}^{\left(k\right)}\left(\frac{1}{\lambda}\right)_{n-m+1}\lambda^{n-m+1}\right.\\
&\left.-\lambda\sum_{m=0}^{n-1}\binom{n}{m}\left(\frac{1}{\lambda}\right)_{n-m}\lambda^{n-m}m\beta_{m,\lambda}^{\left(k\right)}\right\} .
\end{align*}

\end{thm}
Note that
\[
B_{n}^{\left(k\right)}=\lim_{\lambda\rightarrow0}\beta_{n,\lambda}^{\left(k\right)}=\frac{1}{n+1}\left\{ B_{n}^{\left(k-1\right)}-\sum_{m=1}^{n-1}\binom{n}{m-1}B_{m}^{\left(k\right)}\right\} .
\]

Now, we observe that
\begin{align}
 & \relphantom{=}\sum_{n=0}^{\infty}\left(1-\left(1+\lambda t\right)^{-\frac{1}{\lambda}}\right)^{n}\left(n+1\right)^{-k}\label{eq:32}\\
 & =\sum_{n=1}^{\infty}\frac{\left(1-\left(1+\lambda t\right)^{-\frac{1}{\lambda}}\right)^{n}}{n^{k}}\frac{1}{1-\left(1+\lambda t\right)^{-\frac{1}{\lambda}}}\nonumber \\
 & =\frac{1}{1-\left(1+\lambda t\right)^{-\frac{1}{\lambda}}}\li_{k}\left(1-\left(1+\lambda t\right)^{-\frac{1}{\lambda}}\right)\nonumber \\
 & =\sum_{n=0}^{\infty}\beta_{n,\lambda}^{\left(k\right)}\frac{t^{n}}{n!}.\nonumber
\end{align}

By (\ref{eq:32}), we get
\begin{align}
 & \relphantom{=}\sum_{k=0}^{\infty}\left(\sum_{n=0}^{\infty}\beta_{n,\lambda}^{\left(-k\right)}\frac{x^{n}}{n!}\right)\frac{y^{k}}{k!}\label{eq:33}\\
 & =\sum_{k=0}^{\infty}\left(\sum_{m=0}^{\infty}\left(1-\left(1+\lambda t\right)^{-\frac{1}{\lambda}}\right)^{m}\left(m+1\right)^{k}\right)\frac{y^{k}}{k!}\nonumber \\
 & =\sum_{m=0}^{\infty}\left(1-\left(1+\lambda x\right)^{-\frac{1}{\lambda}}\right)^{m}\sum_{k=0}^{\infty}\left(m+1\right)^{k}\frac{y^{k}}{k!}\nonumber \\
 & =\sum_{m=0}^{\infty}\left(1-\left(1+\lambda x\right)^{-\frac{1}{\lambda}}\right)^{m}e^{\left(m+1\right)y}\nonumber \\
 & =\sum_{j=0}^{\infty}\left(-1\right)^{j}\left(e^{-\frac{1}{\lambda}\log\left(1+\lambda x\right)}-1\right)^{j}e^{\left(j+1\right)y}\nonumber \\
 & =\sum_{j=0}^{\infty}\left(-1\right)^{j}j!\sum_{m=j}^{\infty}S_{2}\left(m,j\right)\left(-1\right)^{m}\lambda^{-m}\frac{\left(\log\left(1+\lambda x\right)\right)^{m}}{m!}e^{\left(j+1\right)y}\nonumber \\
 & =\sum_{m=0}^{\infty}\sum_{j=0}^{m}\left(-1\right)^{j+m}j!S_{2}\left(m,j\right)\lambda^{-m}\sum_{n=m}^{\infty}S_{1}\left(n,m\right)\lambda^{n}\frac{x^{n}}{n!}e^{\left(j+1\right)y}\nonumber \\
 & =\sum_{n=0}^{\infty}\left(\sum_{m=0}^{n}\sum_{j=0}^{m}\left(-1\right)^{j+m}j!S_{2}\left(m,j\right)\lambda^{n-m}S_{1}\left(n,m\right)e^{\left(j+1\right)y}\right)\frac{x^{n}}{n!}\nonumber \\
 & =\sum_{k=0}^{\infty}\sum_{n=0}^{\infty}\left(\sum_{m=0}^{n}\sum_{j=0}^{m}\left(-1\right)^{j+m}j!\lambda^{n-m}S_{2}\left(m,j\right)S_{1}\left(n,m\right)\left(j+1\right)^{k}\right)\frac{x^{n}}{n!}\frac{y^{k}}{k!}.\nonumber
\end{align}

Therefore, by (\ref{eq:33}), we obtain the following theorem.
\begin{thm}
\label{thm:6} For $k\in\mathbb{Z}$ and $n\ge0$, we have
\[
\beta_{n,\lambda}^{\left(-k\right)}=\sum_{m=0}^{n}\sum_{j=0}^{m}\left(-1\right)^{j+m}j!\lambda^{n-m}\left(j+1\right)^{k}S_{2}\left(m,j\right)S_{1}\left(n,m\right).
\]

\end{thm}
Note that
\[
B_{n}^{\left(-k\right)}=\lim_{\lambda\rightarrow0}\beta_{n,\lambda}^{\left(-k\right)}=\sum_{j=0}^{n}\left(-1\right)^{j+n}j!\left(j+1\right)^{k}S_{2}\left(n,j\right).
\]

From Theorem \ref{thm:1}, we have
\begin{align*}
\frac{d}{dx}\beta_{n,\lambda}^{\left(k\right)}\left(x\right) & =\sum_{l=0}^{n}\binom{n}{l}\beta_{l,\lambda}^{\left(k\right)}\frac{d}{dx}\left(\prod_{i=0}^{n-l-1}\left(x-i\lambda\right)\right)\\
 & =\sum_{l=0}^{n}\binom{n}{l}\beta_{l,\lambda}^{\left(k\right)}\sum_{j=0}^{n-l-1}\frac{1}{\left(x-j\lambda\right)}\prod_{i=0}^{n-l-1}\left(x-i\lambda \right).
\end{align*}

The generalized falling factorial $\left(x\mid\lambda\right)_{n}$ is given by
\begin{equation}
\left(x\mid\lambda\right)_{n}=x\left(x-\lambda\right)\left(x-2\lambda\right)\cdots\left(x-\left(n-1\right)\lambda\right),\quad\left(n\ge0\right).\label{eq:34}
\end{equation}

As is well known, the Bernoulli numbers of the second kind are defined
by the generating funcdtion
\begin{equation}
\frac{t}{\log\left(1+t\right)}=\sum_{n=0}^{\infty}b_{n}\frac{t^{n}}{n!},\quad\left(\text{see \cite{key-20}}\right).\label{eq:35}
\end{equation}

We observe that

\begin{align}
\int_{0}^{1}\left(1+\lambda t\right)^{\frac{x}{\lambda}}dx & =\sum_{n=0}^{\infty}\lambda^{n}\int_{0}^{1}\left(\frac{x}{\lambda}\right)_{n}dx\frac{t^{n}}{n!}\label{eq:36}\\
 & =\sum_{n=0}^{\infty}\int_{0}^{1}\left(x\mid\lambda\right)_{n}dx\frac{t^{n}}{n!}.\nonumber
\end{align}

On the other hand,
\begin{align}
 & \relphantom{=}\int_{0}^{1}\left(1+\lambda t\right)^{\frac{x}{\lambda}}dx\label{eq:37}\\
 & =\frac{\lambda}{\log\left(1+\lambda t\right)}\left(\left(1+\lambda t\right)^{\frac{1}{\lambda}}-1\right)\nonumber \\
 & =\frac{\lambda t}{\log\left(1+\lambda t\right)}\left(\frac{\left(1+\lambda t\right)^{\frac{1}{\lambda}}-1}{t}\right)\nonumber \\
 & =\left(\sum_{m=0}^{\infty}b_{m}\lambda^{m}\frac{t^{m}}{m!}\right)\left(\sum_{l=0}^{\infty}\frac{\left(1|\lambda \right)_{l+1}}{l+1}\frac{t^{l}}{l!}\right)\nonumber \\
 & =\sum_{n=0}^{\infty}\left(\sum_{l=0}^{n}\frac{\left(1\mid\lambda\right)_{l+1}}{l+1}\lambda^{n-l}b_{n-l}\binom{n}{l}\right)\frac{t^{n}}{n!}.\nonumber
\end{align}

From (\ref{eq:36}) and (\ref{eq:37}), we have
\begin{equation}
\int_{0}^{1}\left(x\mid\lambda\right)_{n}dx=\sum_{l=0}^{n}\binom{n}{l}\lambda^{n-l}b_{n-l}\frac{\left(1\mid\lambda\right)_{l+1}}{l+1},\quad\left(n\ge0\right).\label{eq:38}
\end{equation}

By Theorem \ref{thm:1}, we get

\begin{align*}
\int_{0}^{1}\beta_{n,\lambda}^{\left(k\right)}\left(x\right)dx & =\sum_{l=0}^{n}\binom{n}{l}\beta_{l,\lambda}^{\left(k\right)}\int_{0}^{1}\left(\frac{x}{\lambda}\right)_{n-l}\lambda^{n-l}dx\\
 & =\sum_{l=0}^{n}\binom{n}{l}\beta_{n-l,\lambda}^{\left(k\right)}\int_{0}^{1}\left(x\mid\lambda\right)_{l}dx\\
 & =\sum_{l=0}^{n}\left(\sum_{m=0}^{l}\binom{l}{m}\lambda^{l-m}b_{l-m}\frac{\left(1\mid\lambda\right)_{m+1}}{m+1}\right)\binom{n}{l}\beta_{n-l,\lambda}^{\left(k\right)}.
\end{align*}

\section{Further remarks}

Let $\mathbb{C}$ be complex number field and let $\mathcal{F}$ be
the set of all formal power series in the variable $t$ over $\mathbb{C}$
with
\begin{equation}
\mathcal{F}=\left\{ \left.f\left(t\right)=\sum_{k=0}^{\infty}a_{k}\frac{t^{k}}{k!}\right|a_{k}\in\mathbb{C}\right\} .\label{eq:39}
\end{equation}

Let $\mathbb{P}$ be the algebra of polynomials in a single variable
$x$ over $\mathbb{C}$ and let $\mathbb{P}^{*}$ be the vector space
of all linear functionals on $\mathbb{P}$. The action of linear functional
$L\in\mathbb{P}^{*}$ on a polynomial $p\left(x\right)$ is denoted
by $\acl L{p\left(x\right)},$ and linearly extended as
\[
\acl{cL+c^{\prime}L^{\prime}}{p\left(x\right)}=c\acl L{p\left(x\right)}+c^{\prime}\acl{L^{\prime}}{p\left(x\right)},
\]
where $c,c^{\prime}\in\mathbb{C}$.

For $f\left(t\right)=\sum_{k=0}^{\infty}a_{k}\frac{t^{k}}{k!}$, we
define a linear functional on $\mathbb{P}$ by setting
\begin{equation}
\acl{f\left(t\right)}{x^{n}}=a_{n}\quad\text{for all }n\ge0.\label{eq:40}
\end{equation}
Thus, by (\ref{eq:40}), we get
\begin{equation}
\acl{t^{k}}{x^{n}}=n!\delta_{n,k},\quad\left(n,k\ge0\right),\quad\left(\text{see \cite{key-5,key-15,key-20}}\right).\label{eq:41}
\end{equation}

For $f_{L}\left(t\right)=\sum_{k=0}^{\infty}\acl L{x^{k}}\frac{t^{k}}{k!}$,
by (\ref{eq:41}), we get $\acl{f_{L}\left(t\right)}{x^{n}}=\acl L{x^{n}}.$
In addition, the mapping $L\mapsto f_{L}\left(t\right)$ is a vector
space isomorphism from $\mathbb{P}^{*}$ onto $\mathcal{F}$. Henceforth,
$\mathcal{F}$ denotes both the algebra of the formal power series
in $t$ and the vector space of all linear functionals on $\mathbb{P}$
and so an element $f\left(t\right)$ of $\mathcal{F}$ can be regarded
as both a formal power series and a linear functional. We refer to $\mathcal{F}$
umbral algebra. The umbral calculus is the study of umbral algebra
(see \cite{key-6,key-14,key-20}). The order $o\left(f\left(t\right)\right)$
of the non-zero power series $f\left(t\right)$ is the smallest integer
$k$ for which the coefficient of $t^{k}$ does not vanish.

If $o\left(f\left(t\right)\right)=1$(respectively, $o\left(f\left(t\right)\right)=0$),
then $f\left(t\right)$ is called a delta (respectively, an invertible)
series (see \cite{key-20}). For $o\left(f\left(t\right)\right)=1$
and $o\left(g\left(t\right)\right)=0$, there exists a unique sequence
$s_{n}\left(x\right)$ of polynomials such that $\acl{g\left(t\right)f\left(t\right)^{k}}{s_{n}\left(x\right)}=n!\delta_{n,k}$,
$\left(n,k\ge0\right)$.

The sequence $s_{n}\left(x\right)$ is called the Sheffer sequence
for $\left(g\left(t\right),f\left(t\right)\right)$, and we write
$s_{n}\left(x\right)\sim\left(g\left(t\right),f\left(t\right)\right)$
(see \cite{key-20}).

Let $f\left(t\right)\in\mathcal{F}$ and $p\left(x\right)\in\mathbb{P}$.
Then, by (\ref{eq:41}), we get
\begin{equation}
f\left(t\right)=\sum_{k=0}^{\infty}\acl{f\left(t\right)}{x^{k}}\frac{t^{k}}{k!},\quad p\left(x\right)=\sum_{k=0}^{\infty}\acl{t^{k}}{p\left(x\right)}\frac{x^{k}}{k!}.\label{eq:42}
\end{equation}

From (\ref{eq:42}), we have
\begin{equation}
p^{\left(k\right)}\left(0\right)=\acl{t^{k}}{p\left(x\right)}=\acl 1{p^{\left(k\right)}\left(x\right)},\label{eq:43}
\end{equation}
where $p^{k}\left(x\right)=\frac{d^{k}}{dx^{k}}p\left(x\right)$,
(see \cite{key-11,key-12,key-20}).

By (\ref{eq:43}), we easily get
\begin{equation}
t^{k}p\left(x\right)=p^{k}\left(x\right),\quad e^{yt}p\left(x\right)=p\left(x+y\right),\quad\text{and }\acl{e^{yt}}{p\left(x\right)}=p\left(y\right).\label{eq:44}
\end{equation}
From (\ref{eq:44}), we have
\begin{equation}
\frac{e^{yt}-1}{t}p\left(x\right)=\int_{x}^{x+y}p\left(u\right)du, \quad\acl{e^{yt}-1}{p\left(x\right)}=p\left(y\right)-p\left(0\right).\nonumber
\end{equation}

Let $f\left(t\right)$ be the linear functional such that
\begin{equation}
\acl{f\left(t\right)}{p\left(x\right)}=\int_{0}^{y}p\left(u\right)du,\label{eq:45}
\end{equation}
for all polynomials $p\left(x\right)$. Then it can be determined by (\ref{eq:42})
to be
\begin{align}
f\left(t\right) & =\sum_{k=0}^{\infty}\frac{\acl{f\left(t\right)}{x^{k}}}{k!}t^{k}\label{eq:46}\\
 & =\sum_{k=0}^{\infty}\frac{y^{k+1}}{\left(k+1\right)!}t^{k}\nonumber \\
 & =\frac{1}{t}\left(e^{yt}-1\right).\nonumber
\end{align}

Thus, for $p\left(x\right)\in\mathbb{P}$, we have
\begin{equation}
\acl{\frac{e^{yt}-1}{t}}{p\left(x\right)}=\int_{0}^{y}p\left(u\right)du. \label{eq:47}
\end{equation}

It is known that
\begin{equation}
s_{n}\left(x\right)\sim\left(g\left(t\right),f\left(t\right)\right)\iff\frac{1}{g\left(\overline{f}\left(t\right)\right)}e^{x\overline{f}\left(t\right)}=\sum_{k=0}^{\infty}s_{k}\left(x\right)\frac{t^{k}}{k!},\quad\left(x\in\mathbb{C}\right)\label{eq:48}
\end{equation}
where $\overline{f}\left(t\right)$ is the compositional inverse of
$f\left(t\right)$ such that $f\left(\overline{f}\left(t\right)\right)=\overline{f}\left(f\left(t\right)\right)=t$
(see \cite{key-11,key-20}).

From (\ref{eq:15}), we note that
\begin{equation}
\beta_{n,\lambda}^{\left(k\right)}\left(x\right)\sim\left(\frac{1-e^{-t}}{\li_{k}\left(1-e^{-t}\right)},\frac{1}{\lambda}\left(e^{\lambda t}-1\right)\right).\label{eq:49}
\end{equation}

That is,
\[
\sum_{n=0}^{\infty}\beta_{n,\lambda}^{\left(k\right)}\left(x\right)\frac{t^{n}}{n!}=\frac{\li_{k}\left(1-\left(1+\lambda t\right)^{-\frac{1}{\lambda}}\right)}{1-\left(1+\lambda t\right)^{-\frac{1}{\lambda}}}\left(1+\lambda t\right)^{\frac{x}{\lambda}}.
\]

Thus, by (\ref{eq:49}),
\begin{equation}
\frac{1}{\lambda}\left(e^{\lambda t}-1\right)\beta_{n,\lambda}^{\left(k\right)}\left(x\right)=n\beta_{n-1,\lambda}^{\left(k\right)}\left(x\right).\label{eq:50}
\end{equation}

On the other hand,
\begin{equation}
\left(e^{\lambda t}-1\right)\beta_{n,\lambda}^{\left(k\right)}\left(x\right)=\beta_{n,\lambda}^{\left(k\right)}\left(x+\lambda\right)-\beta_{n,\lambda}^{\left(k\right)}\left(x\right).\label{eq:51}
\end{equation}

Therefore, by (\ref{eq:50}) and (\ref{eq:51}), we obtain the following
theorem.
\begin{thm}
\label{thm:7} For $n\in\mathbb{N}$, we have
\[
\lambda\beta_{n-1,\lambda}^{\left(k\right)}\left(x\right)=\frac{1}{n}\left\{ \beta_{n,\lambda}^{\left(k\right)}\left(x+\lambda\right)-\beta_{n,\lambda}^{\left(k\right)}\left(x\right)\right\} .
\]

\end{thm}
By (\ref{eq:47}), we get
\begin{equation}
\frac{e^{yt}-1}{t}\beta_{n,\lambda}^{\left(k\right)}\left(x\right)=\int_{x}^{x+y}\beta_{n,\lambda}^{\left(k\right)}\left(u\right)du.\label{eq:52}
\end{equation}

From (\ref{eq:52}), we have
\begin{equation}
\acl{\frac{e^{yt}-1}{t}}{\beta_{n,\lambda}^{\left(k\right)}\left(x\right)}=\int_{0}^{y}\beta_{n,\lambda}^{\left(k\right)}\left(u\right)du.\label{eq:53}
\end{equation}

Thus, by (\ref{eq:53}), we get
\begin{align}
 & \relphantom{=}\acl{\frac{e^{t}-1}{t}}{\beta_{n,\lambda}^{\left(k\right)}\left(x\right)}\label{eq:54}\\
 & =\int_{0}^{1}\beta_{n,\lambda}^{\left(k\right)}\left(u\right)du\nonumber \\
 & =\sum_{l=0}^{n}\sum_{m=0}^{l}\binom{l}{m}\binom{n}{l}\lambda^{l-m}b_{l-m}\beta_{n-l,\lambda}^{\left(k\right)}\frac{\left(1\mid\lambda\right)_{m+1}}{m+1}.\nonumber
\end{align}

Therefore, by (\ref{eq:54}), we obtain the following theorem.,
\begin{thm}
\label{thm:8} For $n\ge0$, we have
\[
\acl{\frac{e^{t}-1}{t}}{\beta_{n,\lambda}^{\left(k\right)}\left(x\right)}=\sum_{l=0}^{n}\sum_{m=0}^{l}\binom{l}{m}\binom{n}{l}\lambda^{l-m}b_{l-m}\beta_{n-l,\lambda}^{\left(k\right)}\frac{\left(1\mid\lambda\right)_{m+1}}{m+1}.
\]

\end{thm}
Note that
\begin{align*}
 & \relphantom{=}\acl{\frac{e^{t}-1}{t}}{B_{n}^{\left(k\right)}\left(x\right)}\\
 & =\lim_{\lambda\rightarrow0}\acl{\frac{e^{t}-1}{t}}{\beta_{n,\lambda}^{\left(k\right)}\left(x\right)}\\
 & =\lim_{\lambda\rightarrow0}\int_{0}^{1}\beta_{n,\lambda}^{\left(k\right)}\left(u\right)du\\
 & =\sum_{l=0}^{n}\binom{n}{l}B_{n-l}^{\left(k\right)}\frac{1}{l+1}
\end{align*}

Let
\[
\mathbb{P}_{n}=\left\{ \left.p\left(x\right)\in\mathbb{C}\left[x\right]\right|\deg p\left(x\right)\le n\right\} ,\quad\left(n\ge0\right).
\]

For $p\left(x\right)\in\mathbb{P}_{n}$ with $p\left(x\right)=\sum_{m=0}^{n}a_{m}\beta_{m,\lambda}^{\left(k\right)}\left(x\right)$,
we have
\begin{align}
 & \relphantom{=}\acl{\frac{1-e^{-t}}{\li_{k}\left(1-e^{-t}\right)}\left(\frac{1}{\lambda}\left(e^{\lambda t}-1\right)\right)^{m}}{p\left(x\right)}\label{eq:55}\\
 & =\sum_{l=0}^{n}a_{l}\acl{\frac{1-e^{-t}}{\li_{k}\left(1-e^{-t}\right)}\left(\frac{1}{\lambda}\left(e^{\lambda t}-1\right)\right)^{m}}{\beta_{l,\lambda}^{\left(k\right)}\left(x\right)}\nonumber
\end{align}

From (\ref{eq:49}), we note that
\begin{equation}
\acl{\frac{1-e^{-t}}{\li_{k}\left(1-e^{-t}\right)}\left(\frac{1}{\lambda}\left(e^{\lambda t}-1\right)\right)^{m}}{\beta_{l,\lambda}^{\left(k\right)}\left(x\right)}=l!\delta_{l,m}.\label{eq:56}
\end{equation}

By (\ref{eq:55}) and (\ref{eq:56}), we get
\begin{equation}
a_{m}=\frac{1}{m!}\acl{\frac{1-e^{-t}}{\li_{k}\left(1-e^{-t}\right)}\left(\frac{1}{\lambda}\left(e^{\lambda t}-1\right)\right)^{m}}{p\left(x\right)},\quad\left(m\ge0\right).\label{eq:57}
\end{equation}

Therefore, by (\ref{eq:57}), we obtain the following theorem.
\begin{thm}
\label{thm:9} For $p\left(x\right)\in\mathbb{P}_{n}$, we have
\[
p\left(x\right)=\sum_{m=0}^{n}a_{m}\beta_{m,\lambda}^{\left(k\right)}\left(x\right),
\]
where
\[
a_{m}=\frac{1}{m!}\acl{\frac{1-e^{-t}}{\li_{k}\left(1-e^{-t}\right)}\left(\frac{1}{\lambda}\left(e^{\lambda t}-1\right)\right)^{m}}{p\left(x\right)}.
\]

\end{thm}
For example, let us take $p\left(x\right)=B_{n}^{\left(k\right)}\left(x\right)$
$\left(n\ge0\right)$. Then, by Theorem \ref{thm:9}, we have
\begin{equation}
B_{n}^{\left(k\right)}\left(x\right)=\sum_{m=0}^{n}a_{m}\beta_{m,\lambda}^{\left(k\right)}\left(x\right),\label{eq:58}
\end{equation}
where
\begin{align}
a_{m} & =\frac{1}{m!}\acl{\frac{1-e^{-t}}{\li_{k}\left(1-e^{-t}\right)}\left(\frac{1}{\lambda}\left(e^{\lambda t}-1\right)\right)^{m}}{B_{n}^{\left(k\right)}\left(x\right)}\label{eq:59}\\
 & =\frac{1}{m!}\acl{\frac{1-e^{-t}}{\li_{k}\left(1-e^{-t}\right)}\left(\frac{1}{\lambda}\left(e^{\lambda t}-1\right)\right)^{m}}{\frac{\li_{k}\left(1-e^{-t}\right)}{1-e^{-t}}x^{n}}\nonumber \\
 & =\frac{\lambda^{-m}}{m!}\acl{\left(e^{\lambda t}-1\right)^{m}}{x^{n}}=\lambda^{-m}\sum_{l=m}^{\infty}S_{2}\left(l,m\right)\frac{\lambda^{l}}{l!}\acl{t^{l}}{x^{n}}\nonumber \\
 & =\lambda^{n-m}S_{2}\left(n,m\right).\nonumber
\end{align}

From (\ref{eq:58}) and (\ref{eq:59}), we have
\begin{equation}
B_{n}^{\left(k\right)}\left(x\right)=\sum_{m=0}^{n}\lambda^{n-m}S_{2}\left(n,m\right)\beta_{m,\lambda}^{\left(k\right)}\left(x\right).\label{eq:60}
\end{equation}

\bibliographystyle{amsplain}

\begin{thebibliography}{10}

\bibitem{key-1}
S.~Araci, M.~Acikgoz, and A.~Kilicman, \emph{Extended $p$-adic $q$-invariant
  integrals on $\mathbb{Z}_p$ associated with applications of umbral calculus},
  Adv. Difference Equ. (2013), 2013:96, 14 pp.

\bibitem{key-2}
A.~Bayad, and T.~Kim, \emph{Results on values of Barnes polynomials}, Rocky Mountain J. Math.
  \textbf{43} (2013), 1857–-1869.

\bibitem{key-18}
J.~Broedel, O.~Schlotterer and S. Stieberger,  \emph{Polylogarithms, multiple zeta values and superstring amplitudes},
 Fortschr. Phys. \textbf{61} (2013),no.~9, 812--870.

\bibitem{key-3}
L.~Carlitz, \emph{Degenerate {S}tirling, {B}ernoulli and {E}ulerian numbers},
  Utilitas Math. \textbf{15} (1979), 51--88.

\bibitem{key-4}
D.~Ding and J.~Yang, \emph{Some identities related to the {A}postol-{E}uler and
  {A}postol-{B}ernoulli polynomials}, Adv. Stud. Contemp. Math. (Kyungshang)
  \textbf{20} (2010), no.~1, 7--21.

\bibitem{key-5}
H.~Gangl, \emph{Herbert Functional equations and ladders for polylogarithms}, Commun. Number Theory Phys. \textbf{7} (2013),
  no.~3, 397--410.

\bibitem{key-7}
S.~Gaboury, R.~Tremblay and B.-J. Fugere, \emph{Some explicit formulas for
  certain new classes of {B}ernoulli, {E}uler and {G}enocchi polynomials},
  Proc. Jangjeon Math. Soc. \textbf{17} (2014), no.~1, 97--104.

\bibitem{key-8}
Y.~He and W.~Zhang, \emph{A convolution formula for {B}ernoulli polynomials},
  Ars Combin. \textbf{108} (2013), 97--104.


\bibitem{key-6}
H.~Jolany, M.~Aliabadi, R.B.~Corcino and M.R.~Darafsheh,  \emph{A note on multi poly-Euler numbers and
Bernoulli polynomials},General Mathematics \textbf{20} (2012), no.~2--3, 122--134.

\bibitem{key-16}
M.~Kaneko, \emph{ poly-{B}ernoulli numbers}, J Th{\'e}or. Nombres
  Bordeaux \textbf{9} (1997), no.~1, 221--228.

\bibitem{key-10}
D.~S. Kim and T.~Kim, \emph{A note on poly-{B}ernoulli and higher-order
  poly-{B}ernoulli polynomials}, Russ. J. Math. Phys. \textbf{22} (2015),
  no.~1, 26--33.

\bibitem{key-11}
\bysame, \emph{Higher-order {C}auchy of the first kind and
  poly-{C}auchy of the first kind mixed type polynomials}, Adv. Stud. Contemp.
  Math. (Kyungshang) \textbf{23} (2013), no.~4, 621--636.

\bibitem{key-9}
\bysame, \emph{Higher-order {F}robenius-{E}uler and poly-{B}ernoulli mixed-type
  polynomials}, Adv. Difference Equ. (2013), 2013:251, 13 pp.

\bibitem{key-13}
\bysame, \emph{Daehee polynomials with $q$-parameter}, Adv. Studies Theor.
  Phys. \textbf{8} (2014), no.~13, 561--569.

\bibitem{key-12}
D.~S. Kim, T.~Kim, T.~Mansour, and D.~V. Dolgy, \emph{On poly-{B}ernoulli
  polynomials of the second kind with umbral calculus viewpoint}, Adv.
  Difference Equ. (2015), 2015:27, 13 pp.

\bibitem{key-14}
T.~Kim, \emph{Barnes' type multiple degenerate {B}ernoulli and {E}uler
  polynomials}, Appl. Math. Comput. (2015), 556--564.

\bibitem{key-15}
T.~Kim and T.~Mansour, \emph{Umbral calculus associated with {F}robenius-type
  {E}ulerian polynomials}, Russ. J. Math. Phys. \textbf{21} (2014), no.~4,
  484--493.

\bibitem{key-17}
Q.-M. Luo, \emph{Some recursion fformula and relations for {B}ernoulli numbers
  and {E}uler numbers of higher order}, Adv. Stud. Contemp. Math. (Kyungshang)
  \textbf{10} (2005), no.~1, 63--70.

\bibitem{key-19}
F.~Qi, \emph{An integral representation, complete monotonicity, and
  inequalities of {C}auchy numbers of the second kind}, J. Number Theory
  \textbf{144} (2014), 244--255.

\bibitem{key-20}
S.~Roman, \emph{The umbral calculus}, Pure and Applied Mathematics, vol. 111,
  Academic Press, Inc. [Harcourt Brace Jovanovich, Publishers], New York, 1984.
  \MR{741185 (87c:05015)}

\bibitem{key-21}
E.~Sen, \emph{Theokind on {A}postol-{E}uler polynomials of higher order arising
  from {E}uler basis}, Adv. Stud. Contemp. Math. (Kyungshang) \textbf{23}
  (2013), no.~2, 337--345.

\end{thebibliography}
\providecommand{\bysame}{\leavevmode\hbox to3em{\hrulefill}\thinspace}
\providecommand{\MR}{\relax\ifhmode\unskip\space\fi MR }
\providecommand{\MRhref}[2]{%
  \href{http://www.ams.org/mathscinet-getitem?mr=#1}{#2}
}
\providecommand{\href}[2]{#2}

\end{document}